\documentclass[11pt]{article}
\usepackage{amssymb,amsmath,amscd} 
\usepackage{osid}
\def\mysavedown#1{\edef\mysubs{\mysubs#1}}
\def\mysaveup#1{\edef\mysups{\mysups#1}}
\def\mydown#1{{\mytensor}_{\vphantom{\mysubs}#1}}
\def\myup#1{{\mytensor}^{\vphantom{\mysups}#1}}
\def\tensor#1#2{
  #1
  \def\mytensor{\vphantom{#1}}
  \def\mysubs{\relax}
  \def\mysups{\relax}
  \let\down=\mysavedown
  \let\up=\mysaveup
  #2
  \let\down=\mydown
  \let\up=\myup
  #2
  }

\newcommand{\bbC}{{\mathbb C}} 
\newcommand{\bbR}{{\mathbb R}}

\newcommand{\bbN}{{\mathbb N}} 
\newcommand{\bbK}{{\mathbb K}}

\newcommand{\calH}{{\mathcal H}}

\newcommand{\calP}{{\mathcal P}} 
\newcommand{\calS}{{\mathcal S}} 
\newcommand{\calT}{{\mathcal T}} 
\newcommand{\calV}{{\mathcal V}} 
 
\newcommand{\id}{\mathrm{id}} 
\newcommand{\frakl}{{\mathfrak l}} 
\newcommand{\frakr}{{\mathfrak r}} 

\newcommand{\fR}{{\mathfrak R}}
 
\newcommand{\bpi}{\begin{picture}} 
\newcommand{\epi}{\end{picture}}

\title{Methods for the construction of generators of algebraic curvature tensors} 
\author{Bernd Fiedler\thanks{\hspace*{0.2cm}1991 Mathematics Subject
Classification: 53B20, 15A72, 05E10, 16D60, 05-04.} \\
{\footnotesize\it Eichelbaumstr. 13, D-04249
Leipzig, Germany. E-mail: bfiedler@fiemath.de}\\[2ex]
{\footnotesize\it Dedicated to the memory of Professor Brian G.
Wybourne.}}

\begin{document} 

\maketitle 
	  
\begin{abstract}
We demonstrate the use of several tools from Algebraic Combinatorics such as Young tableaux, symmetry operators, the Littlewood-Richardson rule and discrete Fourier transforms of symmetric groups in investigations of algebraic curvature tensors.
\end{abstract}

In \cite{fie03b,fie03c,fie04a} we constructed and investigated generators of algebraic curvature tensors and algebraic covariant derivative curvature tensors. These investigations followed the example of the paper \cite{full4} by S.A. Fulling, R.C. King, B.G.Wybourne and C.J. Cummins and applied tools from Algebraic Combinatorics such as Young tableaux, symmetry operators (in particular Young symmetrizers), the Littlewood-Richardson rule, but also discrete Fourier transforms of symmetric groups. The present paper is a short summary of \cite{fie03b,fie03c,fie04a} in which we want to demonstrate the use of these methods.

\section{The problem}
Let $V$ be a finite dimensional $\bbK$-vector space, $\bbK = \bbR, \bbC$, and let $\calT_r V$ denote the $\bbK$-vector space of covariant tensors of order $r$ over $V$.
\begin{definition}{Definition}
{\it Algebraic curvature tensors} $\fR\in\calT_4 V$ and {\it algebraic covariant derivative curvature tensors} $\fR'\in\calT_5 V$ are tensors of order 4 or 5 whose coordinates satisfy
\begin{eqnarray*}
{\fR}_{ijkl}\;=\;- {\fR}_{jikl}\;=\;{\fR}_{klij} & \;\;\; & {\fR'}_{ijklm}\;=\;- {\fR'}_{jiklm}\;=\;{\fR'}_{klijm}\\
{\fR}_{ijkl} + {\fR}_{iklj}+ {\fR}_{iljk}\;=\;0 & \;\;\; & {\fR'}_{ijklm} + {\fR'}_{ikljm}+ {\fR'}_{iljkm}\;=\;0\\
 & \;\;\; & {\fR'}_{ijklm} + {\fR'}_{ijlmk} + {\fR'}_{ijmkl}\;=\;0\,.
\end{eqnarray*}
They are tensors which possess the same symmetry properties as the Riemannian curvature tensor $R_{ijkl}$ and its covariant derivative $R_{ijkl;m}$ of a Levi-Civita connection $\nabla$.
\end{definition}
The vector space of algebraic curvature tensors $\fR\in\calT_4 V$
is spanned by each of the following sets of tensors (P. Gilkey \cite[pp.41-44]{gilkey5}, B. Fiedler\footnote{\cite{fie20} uses also tools of Algebraic Combinatorics, in particular plethysms.} \cite{fie20})
\begin{eqnarray} \label{equ1.1}%
\hspace*{1cm}\gamma (S)_{i j k l} & := &
S_{i l} S_{j k} - S_{i k} S_{j l}\,,\hspace{75pt}S\;{\rm symmetric} \\
\alpha (A)_{i j k l} & := &
2\,A_{i j} A_{k l} + A_{i k} A_{j l}
-  A_{i l} A_{j k}\,,\;\;\;\;\;A\;{\rm skew-symmetric} \,. \label{equ1.2}%
\end{eqnarray}
The vector space of algebraic covariant derivative curvature tensors $\fR'\in\calT_5 V$ is spanned by the following set of tensors (P. Gilkey \cite[p.236]{gilkey5}, B. Fiedler \cite{fie03b})
\begin{eqnarray} \label{equ1.3}%
\hspace*{-4mm}\hat{\gamma} (S,\hat{S})_{i j k l s} & := &
S_{il}{\hat{S}}_{jks} - S_{jl}{\hat{S}}_{iks} + S_{jk}{\hat{S}}_{ils} - S_{ik}{\hat{S}}_{jls}\,,\hspace{0.5cm}S\,,\,\hat{S}\;{\rm symmetric}.
\end{eqnarray}
\begin{definition}{Problem} \label{problem1.1}%
In the present paper we search for generators of algebraic curvature tensors $\fR$ or algebraic covariant derivative curvature tensors $\fR'$ which can be formed by a suitable symmetry operator from the following types of tensors
\begin{eqnarray} \label{equ1.4}%
\fR\,:\;\;\;U\otimes w\hspace{0.2cm} &\;\;\;,\;\;\;& U\in\calT_3 V\;,\;w\in\calT_1 V\,,\\
\fR'\,:\;\;\;U\otimes W &\;\;\;,\;\;\;& U\in\calT_3 V\;,\;W\in\calT_2 V\,, \label{equ1.5}%
\end{eqnarray}
where $W$ and $U$ belong to symmetry classes of $\calT_2 V$ and $\calT_3 V$ which are defined by minimal right ideals $\frakr\subset\bbK[\calS_2]$ and $\hat{\frakr}\subset\bbK[\calS_3]$, respectively.
\end{definition}
We use Boerner's definition \cite[p.127]{boerner} of symmetry classes of tensors. An element $a = \sum_{p \in {\calS}_r} a(p)p\in\bbK [{\calS}_r]$ of the {\it group ring} of the {\it symmetric group} $\calS_r$ can be considered a {\it symmetry operator} for covariant tensors $T\in\calT_r V$.
The action of $a$ on $T$ is defined by:
\begin{eqnarray}
(a T)(v_1 , \ldots , v_r) & := & \sum_{p \in {\calS}_r} a(p)\,
T(v_{p(1)}, \ldots , v_{p(r)}) \;\;\;\;\;,\;\;\;\;\;
v_i \in V \,.
\end{eqnarray}
\begin{definition}{Definition}
Let $\frakr \subseteq \bbK [{\calS}_r]$ be a right ideal of $\bbK [{\calS}_r]$ for which an $a \in \frakr$ and a $T \in {\calT}_r V$ exist such that $aT \not= 0$. Then the tensor set
\begin{eqnarray}
{\calT}_{\frakr} & := & \{ a T \;|\; a \in \frakr \;,\;
T \in {\calT}_r V \}
\end{eqnarray}
is called the {\itshape symmetry class} of tensors defined by $\frakr$.
\end{definition}
Boerner \cite[p.127]{boerner} showed: If $e\in\bbK[\calS_r]$ is a generating idempotent of $\frakr$, i.e. $\frakr = e\cdot\bbK[\calS_r]$, then it holds
\begin{center}
\fbox{$T\in\calT_r V$ belongs to $\calT_{\frakr}\;\;\;\Leftrightarrow\;\;\;eT\;=\;T$\,.
}
\end{center}
\section{Young symmetrizers}
Young symmetrizers are important symmetry operators. In particular the symmetries of the Riemann tensor $R$ and its covariant derivatives are characterized by a Young symmetrizer. First we define Young tableaux.

A {\itshape Young tableau} $t$ of $r\in\bbN$ is an arrangement of $r$ boxes such that
\begin{enumerate}
\item{the numbers ${\lambda}_i$ of boxes in the rows $i = 1 , \ldots , l$ form a decreasing sequence
${\lambda}_1 \ge {\lambda}_2 \ge \ldots \ge {\lambda}_l > 0$ with
${\lambda}_1 + \ldots + {\lambda}_l = r$,}
\item{the boxes are fulfilled by the numbers $1, 2, \ldots , r$ in any order.}
\end{enumerate}
For instance, the following graphics shows a Young tableau of $r = 16$.
{\small
\[\left.
\begin{array}{cc|c|c|c|c|c|c}
\cline{3-7}
{\lambda}_1 = 5 & \;\;\; & 11 & 2 & 5 & 4 & 12 & \\
\cline{3-7}
{\lambda}_2 = 4 & \;\;\; & 9 & 6 & 16 & 15 & \multicolumn{2}{c}{\;\;\;} \\
\cline{3-6}
{\lambda}_3 = 4 & \;\;\; & 8 & 14 & 1 & 7 & \multicolumn{2}{c}{\;\;\;} \\
\cline{3-6}
{\lambda}_4 = 2 & \;\;\; & 13 & 3 & \multicolumn{4}{c}{\hspace{2cm}} \\
\cline{3-4}
{\lambda}_5 = 1 & \;\;\; & 10 & \multicolumn{4}{c}{\hspace{2cm}} \\
\cline{3-3}
\end{array}\right\}\;=\;t\,.
\]
}
Obviously, the unfilled arrangement of boxes, the {\itshape Young frame}, is characterized by a partition
$\lambda = ({\lambda}_1 , \ldots , {\lambda}_l) \vdash r$ of $r$.

If a Young tableau $t$ of a partition $\lambda \vdash r$ is given, then the {\itshape Young symmetrizer} $y_t$ of $t$ is defined by
\begin{eqnarray}
y_t & := & \sum_{p \in {\calH}_t} \sum_{q \in {\calV}_t} \mathrm{sign}(q)\, p \circ q
\end{eqnarray}
where ${\calH}_t$, ${\calV}_t$ are the groups of the {\itshape horizontal} or
{\itshape vertical permutations} of $t$ which only permute numbers within rows or columns of $t$, respectively. The Young symmetrizers of $\bbK [{\calS}_r]$ are essentially idempotent and define decompositions
\begin{eqnarray}
\bbK [{\calS}_r] \;=\;
\bigoplus_{\lambda \vdash r} \bigoplus_{t \in {\calS\calT}_{\lambda}}
\bbK [{\calS}_r]\cdot y_t
& \;\;,\;\; &
\bbK [{\calS}_r] \;=\;
\bigoplus_{\lambda \vdash r} \bigoplus_{t \in {\calS\calT}_{\lambda}}
y_t \cdot \bbK [{\calS}_r] \label{equ2.2}
\end{eqnarray}
of $\bbK [{\calS}_r]$ into minimal left or right ideals $\bbK [{\calS}_r]\cdot y_t$, $y_t \cdot \bbK [{\calS}_r]$. In \eqref{equ2.2}, the symbol ${\calS\calT}_{\lambda}$ denotes the set of all standard tableaux of the partition $\lambda$. Standard tableaux are Young tableaux in which the entries of every row and every column form an increasing number sequence.\footnote{About Young symmetrizers and
Young tableaux see for instance
\cite{boerner,full4,fulton,jameskerb,kerber,littlew1,mcdonald}.}

The inner sums of (\ref{equ2.2}) are minimal two-sided ideals
\begin{eqnarray}
{\mathfrak a}_{\lambda} & := &
\bigoplus_{t \in {\calS\calT}_{\lambda}}
\bbK [{\calS}_r]\cdot y_t
\;=\;
\bigoplus_{t \in {\calS\calT}_{\lambda}}
y_t \cdot \bbK [{\calS}_r]
\end{eqnarray}
of $\bbK [{\calS}_r]$.
The set of all Young symmetrizers $y_t$ which lie in ${\mathfrak a}_{\lambda}$ is equal to the set of all $y_t$ whose tableau $t$ has the frame $\lambda\vdash r$.
Furthermore two minimal left ideals $\frakl_1, \frakl_2\subseteq\bbK[\calS_r]$ or two minimal right ideals $\frakr_1, \frakr_2\subseteq\bbK[\calS_r]$ are \textit{equivalent} iff they lie in the same ideal ${\mathfrak a}_{\lambda}$.
Now we say that a symmetry class $\calT_{\frakr}$ \textit{belongs to} $\lambda\vdash r$ iff
$\frakr\subseteq{\mathfrak a}_{\lambda}$.

S.A. Fulling, R.C. King, B.G.Wybourne and C.J. Cummins showed in \cite{full4} that the symmetry classes of the Riemann tensor $R$ and its {\itshape symmetrized\footnote{$(\,\ldots\,)$ denotes the symmetrization with respect to the indices $s_1, \ldots , s_u$.} covariant derivatives}
\begin{eqnarray}
\left({\nabla}^{(u)}R\right)_{i j k l s_1 \ldots s_u} & := & {\nabla}_{(s_1} {\nabla}_{s_2} \ldots {\nabla}_{s_u)} R_{i j k l}\;=\;R_{i j k l\,;\,(s_1 \ldots s_u)}
\end{eqnarray}
are generated by special Young symmetrizers\footnote{A proof of this result of \cite{full4} can be found in \cite[Sec.6]{fie12}, too. See also \cite{fie03b} for more details.}. 
\begin{theorem}{Proposition} {\rm\bf (Fulling, King, Wybourne, Cummins)} \label{thm2.3}\\
Let $\nabla$ be the Levi-Civita connection of a pseudo-Riemannian metric $g$.
For $u \ge 0$ the symmetrized covariant derivatives
${\nabla}^{(u)} R$ fulfil
\begin{eqnarray} \label{e12}%
e_t^{\ast} {\nabla}^{(u)} R & = & {\nabla}^{(u)} R
\end{eqnarray}
where $e_t := y_t (u+1)/(2\cdot (u+3)!)$ is an idempotent which is formed from the Young symmetrizer $y_t$ of the standard tableau
\begin{eqnarray} \label{equ2.12}%
t & = &
\begin{array}{|c|c|c|cc|c|}
\hline
1 & 3 & 5 & \ldots & \ldots & (u+4) \\
\hline
2 & 4 & \multicolumn{4}{l}{\hspace{3cm}} \\
\cline{1-2}
\end{array} \,.
\end{eqnarray}
The '$\ast$' in {\rm (\ref{e12})} is the mapping $\ast: a = \sum_{p \in {\calS}_r} a(p)\,p \;\mapsto\; a^{\ast} :=
\sum_{p \in {\calS}_r} a(p)\,p^{-1}$.
\end{theorem}
We see from Proposition \ref{thm2.3} that the tensor fields ${\nabla}^{(u)}R$ belong to the symmetry class which is defined by the symmetrizer $y_t^{\ast}$ of (\ref{equ2.12}), more precisely, by the right ideal
$\frakr = y_t^{\ast}\cdot\bbK [\calS_{u+4}]$. In the special case of algebraic tensors $\fR$, $\fR'$ we have the following corollary (see \cite{fie03b}):
\begin{theorem}{Corollary} \label{cor5}%
Let us denote by $t$ and $t'$ the standard tableaux
\begin{eqnarray} \label{equ2.3}%
t\;=\;
\begin{array}{|c|c|}
\hline
1 & 3 \\
\hline
2 & 4 \\
\hline
\end{array}
& \;\;\;,\;\;\;
t'\;=\;
\begin{array}{|c|c|c|c}
\cline{1-3}
1 & 3 & 5 & \\
\cline{1-3}
2 & 4 &\multicolumn{2}{c}{\;\;\;} \\
\cline{1-2}
\end{array}
\,.
\end{eqnarray}
Then a tensor $T \in {\calT}_4 V$ or $\tilde{T} \in {\calT}_5 V$ is an algebraic curvature tensor or an algebraic covariant derivative curvature tensor iff these tensors satisfy
\begin{eqnarray}
y_t^{\ast} T\;=\;12\,T & \;\;\;,\;\;\; &
y_{t'}^{\ast} \tilde{T}\;=\;24\,\tilde{T}\,,
\end{eqnarray}
respectively. Thus the symmetry classes of the algebraic tensors $\fR$, $\fR'$ are defined by the minimal right ideals $y_t^{\ast}\cdot\bbK[\calS_4]$, $y_{t'}^{\ast}\cdot\bbK[\calS_5]$ which belong to the partitions $(2\,2)\vdash 4$, $(3\,2)\vdash 5$.
\end{theorem}

\section{Symmetry classes of $\calT_3 V$ belonging to $\lambda = (2\,1)$}
The group ring $\bbK[\calS_3]$ contains the minimal 2-sided ideals ${\mathfrak a}_{(3)}, {\mathfrak a}_{(2\,1)}, {\mathfrak a}_{(1^3)}$. The 2-sided ideals ${\mathfrak a}_{(3)}, {\mathfrak a}_{(1^3)}\subset\bbK[\calS_3]$ have dimension\footnote{The dimension of a minimal left or right ideal can be calculated from the Young frame belonging to it by the hook length formula (see e.g. \cite{boerner,full4,jameskerb,kerber}.)} 1 and define consequently unique symmetry classes of $\calT_3 V$. The 2-sided ideal ${\mathfrak a}_{(2\,1)}\subset\bbK[\calS_3]$, however, has dimension 4 and contains an infinite set of minimal right ideals $\hat{\frakr}$ (of dimension 2) which lead to an infinite set of possible symmetry classes for the tensor $U\in\calT_3 V$.

We use diskrete Fourier transforms to determine a generating idempotent for every such $\hat{\frakr}$.
\begin{definition}{Definition}
A {\it discrete Fourier transform}\footnote{Three discrete Fourier transforms are known for symmetric groups $\calS_r$: {\it Young's natural representation} \cite[pp.102-108]{boerner}, {\it Young's seminormal representation} \cite{boerner2}, \cite[p.130]{clausbaum1}, \cite[p.76]{kerber} and {\it Young's orthogonal representation} \cite{boerner2,clausbaum1,kerber} of $\calS_r$.}  for $\calS_r$ is an isomorphism
\begin{eqnarray} \label{e16}%
\hspace*{-5mm}D:\,\bbK[\calS_r]\;\rightarrow\;
\bigotimes_{\lambda \vdash r} {\bbK}^{d_{\lambda} \times d_{\lambda}}
& , &
\sum_{p\in\calS_r}\,a(p)\,p\;\mapsto\;
\left(
\begin{array}{cccc}
A_{{\lambda}_1} & & & 0 \\
 & A_{{\lambda}_2} & & \\
 & & \ddots & \\
0 & & & A_{{\lambda}_k} \\
\end{array}
\right)
\end{eqnarray}
according to Wedderburn's theorem which maps the group ring $\bbK [\calS_r]$ onto an outer direct product $\bigotimes_{\lambda \vdash r} {\bbK}^{d_{\lambda} \times d_{\lambda}}$ of full matrix rings ${\bbK}^{d_{\lambda} \times d_{\lambda}}$. 
\end{definition}
In (\ref{e16}) the matrix ring ${\bbK}^{d_{\lambda} \times d_{\lambda}}$ corresponds to the minimal two-sided ideal $\mathfrak{a}_{\lambda}$ of $\bbK[\calS_r]$.
For $\calS_3$ we have a mapping
\begin{eqnarray}
D :\;\;a\;=\;\sum_{p\in\calS_r}\,a(p)\,p & \mapsto &
\left(
\begin{array}{ccc}
A_{(3)} & & 0 \\
 & A_{(2\,1)} & \\
0 & & A_{(1^3)} \\
\end{array}
\right)\,,
\end{eqnarray}
where $A_{(3)}$ and $A_{(1^3)}$ are $1\times 1$-matrices and $A_{(2\,1)}$ is a $2\times 2$-matrix. It holds $a\in{\mathfrak a}_{(2\,1)}$ iff $A_{(3)} = A_{(1^3)} = 0$. In \cite{fie03b} we proved
\begin{theorem}{Proposition}
Every minimal right ideal $\frakr\subset\bbK^{2\times 2}$ is generated by exactly one of the following (primitive) idempotents
\begin{eqnarray} \label{equ3.3}
X_{\infty}\;:=\;\left(
\begin{array}{cc}
0 & 0 \\
0 & 1 \\
\end{array}
\right)
 & \;\;\;,\;\;\; &
X_{\nu}\;:=\;\left(
\begin{array}{cc}
1 & 0 \\
\nu & 0 \\
\end{array}
\right)\;\;\;,\;\;\;\nu\in\bbK\,.
\end{eqnarray}
\end{theorem}
From (\ref{equ3.3}) we obtain the generating idempotents for the right ideals $\hat{\frakr}\subset{\mathfrak a}_{(2\,1)}\subset\bbK[\calS_3]$ by
\begin{eqnarray} \label{e19}%
{\xi}_{\nu} & \!\!\!=\!\!\! &
D^{-1}\left(
\begin{array}{ccc}
0 & &  \\
 & X_{\nu} & \\
 & & 0 \\
\end{array}
\right) = 
\left\{
\begin{array}{l}
\frac{1}{3}\{[1,2,3] - [2,1,3] - [2,3,1] + [3,2,1]\}, 
\hspace{2mm}\nu = \infty\\*[0.3cm]
\frac{1}{3}\{[1,2,3] + \nu [1,3,2] + (1-\nu)[2,1,3]\\
-\nu [2,3,1] + (-1+\nu)[3,1,2] - [3,2,1]\},\hspace{2.5mm}{\rm else}.\\
\end{array}
\right.
\end{eqnarray}
{\bf Note} that the above connection between $\nu$ and ${\xi}_{\nu}$ depends on the concrete discrete Fourier transform D which is used in (\ref{e19}). The above formula (\ref{e19}) was determined by means of the {\it Mathematica} package {\sf PERMS} \cite{fie10} whose discrete Fourier transform $D$ is based on {\it Young's natural representation} of $\calS_r$.

Examples of tensor fields with a $(2\,1)$-symmetry can be constructed from covariant derivatives of symmetric or alternating 2-fields.
\begin{theorem}{Proposition}\footnote{See \cite{fie03a}.}
Let $\nabla$ be a torsion-free covariant derivative on a $C^{\infty}$-man\-ifold $M$, $\dim M\ge 2$. Further let $\psi,\omega\in\calT_2 M$ be differentiable tensor fields of order {\rm 2} which are symmetric or skew-symmetric, respectively. Then for every point $p\in M$ the tensors
\begin{eqnarray}
(\nabla\psi - \mathrm{sym}(\nabla\psi))|_p\;\Rightarrow\;\nu = 0
& \;\;\;,\;\;\; &
(\nabla\omega - \mathrm{alt}(\nabla\omega))|_p\;\Rightarrow\;\nu = 2
\end{eqnarray}
lie in $(2\,1)$-symmetry classes whose generating idempotents ${\xi}_{\nu}$ belong to the above $\nu$-values $0$ or $2$. {\rm (}'{\rm sym}' = symmetrization, '{\rm alt}' = anti-symmetrization.{\rm )}
\end{theorem}
In a future paper we will show that tensor fields with a $(2\,1)$-symmetry also occur in curvature formulas connected with {\it stationary} and {\it static space-times}.

\section{Our main results}
Now we formulate our main results which were proved in \cite{fie03b,fie04a}. The next section will give some ideas of the proofs.
\begin{theorem}{Theorem} \label{thm3.1}%
A solution of Problem {\rm\ref{problem1.1}} can be constructed at most from such products {\rm (\ref{equ1.4})}or {\rm (\ref{equ1.5})} whose factors belong to the following symmetry classes:
\begin{center}
{\rm
\begin{tabular}{|l|c|l|l|}
\hline
product & & partitions of the symm. classes & \\
\hline
$\fR:\,U\otimes w$ & (a) & $U\rightarrow\,(2\,1)$ & \\
\hline
$\fR':\,U\otimes W$ & (a') & $U\rightarrow\,(3)\hspace{0.35cm},\,W\rightarrow\,(2)$ & $U$ and $W$ symmetric \\
 & (b') & $U\rightarrow\,(2\,1)\,,\,W\rightarrow\,(2)$ & $W$ symmetric \\
 & (c') & $U\rightarrow\,(2\,1)\,,\,W\rightarrow\,(1^2)$ & $W$ skew-symmetric \\
\hline
\end{tabular}\,.
}
\end{center}
\end{theorem}
The case (a') of Theorem \ref{thm3.1} is realized in formula (\ref{equ1.3}).
The cases (a), (b') and (c') of Theorem \ref{thm3.1} lead to
\begin{theorem}{Theorem} \label{thm3.3}
The products {\rm (}a{\rm )}, {\rm (}b'{\rm )}, {\rm (}c'{\rm )} lead to generators
\begin{eqnarray*}
 & &
y_t^{\ast}(U\otimes w)
\;\;,\;\;
y_{t'}^{\ast}(U\otimes S)
\;\;,\;\;
y_{t'}^{\ast}(U\otimes A)
\end{eqnarray*}
of the spaces of algebraic tensors $\fR$, $\fR'$ if and only if the generating idempotent ${\xi}_{\nu}$ of the symmetry class of $U$ fulfills\\
\hspace*{7cm}\fbox{$\nu\;\not=\;\frac{1}{2}$}\,.\\*[0.2cm]
Here $t$ and $t'$ are the Young tableaux {\rm (\ref{equ2.3})}.
\end{theorem}

\section{Ideas of the proofs}

\subsection{Use of the Littlewood-Richardson rule (Proof of Theorem \ref{thm3.1})}
When we consider a right ideal $\frakr$ that defines the symmetry class of a product $T_1\otimes T_2$ of tensors of order $r_1$, $r_2$, then we can determine information about the decomposition of $\frakr$ into minimal right ideals by means of {\it Littlewood-Richardson products} (see \cite{littlew1,full4,fie16,fie18,fie03b}). Let $\frakr_1$, $\frakr_2$ be the right ideals defining the symmetry classes of $T_1$, $T_2$. We consider the left ideals $\frakl_i := \frakr_i^{\ast}$ representation spaces of subrepresentations ${\alpha}_i$ of the natural representation of $\calS_{r_i}$. Then we have $\frakr = \frakl^{\ast}$ where the left ideal $\frakl$ is the representation space of the Littlewood-Richardson product ${\alpha}_1{\alpha}_2 := ({\alpha}_1\,\#\,{\alpha}_2)\uparrow\calS_{r_1+r_2}$. ('$\#$' denotes the outer tensor product and '$\uparrow$' the forming of the induced representation.)

For the tensor products (\ref{equ1.4}), (\ref{equ1.5}) we have to calculate the following {\it Littlewood-Richardson products} by means of the {\it Littlewood-Richardson rule}\footnote{See \cite{littlew1,full4,fie16}.}:
{\small
\[
\begin{array}{lrcllrcl}
\fR: & [3][1] & \sim & [4] + [3\,1] & \fR': & [3][2] & \sim & [5] + [3\,2] + [4\,1]\\
 & [2\,1][1] & \sim & [3\,1] + [2^2] + [2\,1^2]\hspace{2.2mm} & & [3][1^2] & \sim & [4\,1] + [3\,1^2]\\
 & [1^3][1] & \sim & [2\,1^2] + [1^4] & & [2\,1][2] & \sim & [3\,2] + [4\,1] + [2^2\,1] + [3\,1^2]\\
 & & & & & [2\,1][1^2] & \sim & [3\,2] + [2^2\,1] + [3\,1^2] + [2\,1^3] \\
 & & & & & [1^3][2] & \sim & [3\,1^2] + [2\,1^3]\\
 & & & & & [1^3][1^2] & \sim & [2^2\,1] + [2\,1^3] + [1^5]\\
\end{array}
\]
}
Only the products $[2\,1][1]$ for $\fR$ and $[3][2]$, $[2\,1][2]$, $[2\,1][1^2]$ for $\fR'$ contain minimal right ideals that belong to the partitions $(2\,2)$ for $\fR$ and $(3\,2)$ for $\fR'$.
We will definitely obtain $y_{t}^{\ast}(U\otimes w) = 0$ or $y_{t'}^{\ast}(U\otimes W) = 0$ if the ideal $\frakr$ or $\frakr'$ of $W\otimes w$ or $U\otimes W$ does not possesses a subideal belonging to $(2\,2)$ or $(3\,2)$, since then
$y_{t}^{\ast}\in\mathfrak{a}_{(2\,2)}$, $y_{t'}^{\ast}\in\mathfrak{a}'_{(3\,2)}$ but $\frakr\cap\mathfrak{a}_{(2\,2)} = 0$, $\frakr'\cap\mathfrak{a}'_{(3\,2)} = 0$.

\subsection{A step of the proof of Theorem \ref{thm3.3}}
Let us consider the example of expressions $y_{t'}^{\ast}(U\otimes S)$ and $y_{t'}^{\ast}(U\otimes A)$.
To treat such expressions we form the following group ring elements of $\bbK[\calS_5]$:
\begin{eqnarray}
{\sigma}_{\nu , \epsilon} & := & y_{t'}^{\ast}\cdot
{\xi}_{\nu}'\cdot {\zeta}_{\epsilon}''\\
{\zeta}_{\epsilon}'' & := & \id + \epsilon\,(4\,5)\;\;\;,\;\;\;
\epsilon\in\{1,-1\}\\
{\xi}_{\nu} & \mapsto & {\xi}_{\nu}'\in\bbK[\calS_5] \label{equ4.14}
\end{eqnarray}
Formula (\ref{equ4.14}) denotes the embedding of the group ring elements ${\xi}_{\nu}\in\bbK[\calS_3]$ into $\bbK[\calS_5]$ which is induced by the mapping
$\calS_3\rightarrow\calS_5\;,\;
[i_1,i_2,i_3]\mapsto [i_1,i_2,i_3,4,5]$. The symmetry operator ${\xi}_{\nu}'\cdot {\zeta}_{\epsilon}''$ maps arbitrary product tensors $T'''\otimes T''$ to products $U\otimes S$ or $U\otimes A$.
Using our Mathematica package {\sf PERMS} \cite{fie10} we verified 
${\sigma}_{\nu , \epsilon}\not= 0\;\Leftrightarrow\;\nu\not=\frac{1}{2}$.
The value $\nu = \frac{1}{2}$ 
has to be excluded since ${\sigma}_{\nu , \epsilon} = 0$ and $\fR' = 0$ in this case.

\section{Shortest formulas for generators of $\fR$ and $\fR'$}
In this section we want to construct generators (\ref{equ1.4}), (\ref{equ1.5}) whose coordinate representation has a minimal number of summands. To this end we determine systems of linear identities which are satisfied by the coordinates of all tensors from the symmetry class of $U$. In \cite[Sec.III.4.1]{fie16} we proved
\begin{theorem}{Proposition}
Let $\frakr\subset\bbK[\calS_r]$ be a $d$-dimensional right ideal that defines a symmetry class $\calT_{\frakr}$ of tensors $T\in\calT_r V$.
If a basis
$\{ h_1 , \ldots , h_d \}$ of the left ideal $\frakl = \frakr^{\ast}$ is known, then every solution $x_p$ of the linear $(d\times r!)$-equation system
\begin{eqnarray} \label{equ5.4}%
\sum_{p \in {\calS}_r} h_i(p)\,x_p & = & 0 \hspace{1cm}(i = 1 , \ldots , d)\,.
\end{eqnarray}
yields the coefficients for a linear identity
\begin{eqnarray} \label{equ6.1}%
\sum_{p\in\calS_r}\,x_p\,T_{i_{p(1)}\ldots i_{p(r)}} & = & 0
\end{eqnarray}
fulfilled by the coordinates of all $T\in\calT_{\frakr}$.
\end{theorem}
For our tensors $U$ the rank of the equation system (\ref{equ5.4}) is equal to $\dim\frakr = 2$. The columns of (\ref{equ5.4}) are numbered by the permutations $p\in\calS_3$.

Now we determine identities (\ref{equ6.1}) for $U$ by the following procedure. We form the system (\ref{equ5.4}) from the idempotent ${\xi}_{\nu}$ by the determination of a basis\footnote{Faster algorithms which determine a basis also for a large $\calS_r$ by means of discrete Fourier transforms were developed in \cite{fie16}.} from the generating set $\{p\cdot {\xi}_{\nu}^{\ast}\;|\;p\in\calS_3\}$ of $\frakl = \frakr^{\ast}$. Then
for every subset $\calP = \{p_1,p_2\}\subset\calS_3$ we check the determinant ${\Delta}_{\calP}$ of the corresponding $(2\times 2)$-submatrix of (\ref{equ5.4}). If ${\Delta}_{\calP}\not= 0$, then we determine identities (\ref{equ6.1}) of the special form
\begin{eqnarray}
0 & = & \sum_{p\in\calP}\,x_p^{(\bar{p})}\,U_{i_{p(1)} i_{p(2)} i_{p(3)}} +
U_{i_{\bar{p}(1)} i_{\bar{p}(2)} i_{\bar{p}(3)}}\hspace{1cm}(\bar{p}\in{\calS_3}\setminus\calP)\,.
\end{eqnarray}
For instance, the set
$\calP = \{[1,2,3]\,,\,[1,3,2]\}$ leads to the determinant
${\Delta}_{\calP}(\nu) = {\textstyle\frac{4}{27}}\,(1 - \nu)(1 + \nu)$
which has the roots ${\nu}_1 = 1$ and ${\nu}_2 = -1$. For $\nu\not\in\{1\,,\,-1\}$ we obtain the identities
\begin{eqnarray} \label{equ5.7}%
 & &
\hspace*{-1cm}\begin{array}{crcrccccccc}
- & \frac{{\nu}^2 - \nu + 1}{{\nu}^2 - 1}\,U_{ijk} & + & \frac{2\nu - 1}{{\nu}^2 - 1}\,U_{ikj}        & +  &         &         &         & U_{kji} & = & 0 \\
  & \frac{{\nu}^2 - 2\nu}{{\nu}^2 - 1}\,U_{ijk} & + & \frac{{\nu}^2 - \nu + 1}{{\nu}^2 - 1}\,U_{ikj} & + &         &         & U_{kij} &         & = & 0 \\
  & \frac{2\nu - 1}{{\nu}^2 - 1}\,U_{ijk}        & - & \frac{{\nu}^2 - \nu + 1}{{\nu}^2 - 1}\,U_{ikj} & + &         & U_{jki} &         &         & = & 0 \\
  & \frac{{\nu}^2 - \nu + 1}{{\nu}^2 - 1}\,U_{ijk} & + & \frac{{\nu}^2 - 2\nu}{{\nu}^2 - 1}\,U_{ikj} & + & U_{jik} &         &         &         & = & 0 \\
\end{array}\;.
\end{eqnarray}
There exist 15 subsets $\calP = \{p_1,p_2\}\subset\calS_3$ and consequently 15 systems (\ref{equ5.7}) for $U$.

Now let us consider the coordinates $\frac{1}{24}(y_{t'}^{\ast}(U\otimes S))_{i j k l r}$ of generators for $\fR'$. If we use (\ref{equ5.7}) to express all coordinates of $U$ by $U_{ijk}$ and $U_{ikj}$ we obtain the following sum of 16 terms.
\[
\begin{array}{rcrc}
\frac{-\left( -1 + 2\,\nu \right) }{24\,\left( -1 + \nu \right) \,\left( 1 + \nu \right) } 
   \,\,\tensor{U}{\down{j}\down{l}\down{r}} 
    \tensor{S}{\down{i}\down{k}} & + &
  \frac{\nu\,\left( -1 + 2\,\nu \right) }
    {24\,\left( -1 + \nu \right) \,\left( 1 + \nu \right) } 
   \,\,\tensor{U}{\down{j}\down{r}\down{l}} 
    \tensor{S}{\down{i}\down{k}} & + \\ 
  \frac{-1 + 2\,\nu}{24\,\left( -1 + \nu \right) \,\left( 1 + \nu \right) } 
   \,\,\tensor{U}{\down{j}\down{k}\down{r}} 
    \tensor{S}{\down{i}\down{l}} & - &
  \frac{\nu\,\left( -1 + 2\,\nu \right) }
    {24\,\left( -1 + \nu \right) \,\left( 1 + \nu \right) } 
   \,\,\tensor{U}{\down{j}\down{r}\down{k}} 
    \tensor{S}{\down{i}\down{l}} & + \\ 
  \frac{-1 + 2\,\nu}{24\,\left( -1 + \nu \right) } 
   \,\,\tensor{U}{\down{j}\down{k}\down{l}} 
    \tensor{S}{\down{i}\down{r}} & - &
  \frac{-1 + 2\,\nu}{24\,\left( -1 + \nu \right) } 
   \,\,\tensor{U}{\down{j}\down{l}\down{k}} 
    \tensor{S}{\down{i}\down{r}} & + \\ 
  \frac{-1 + 2\,\nu}{24\,\left( -1 + \nu \right) \,\left( 1 + \nu \right) } 
   \,\,\tensor{U}{\down{i}\down{l}\down{r}} 
    \tensor{S}{\down{j}\down{k}} & - &
  \frac{\nu\,\left( -1 + 2\,\nu \right) }
    {24\,\left( -1 + \nu \right) \,\left( 1 + \nu \right) } 
   \,\,\tensor{U}{\down{i}\down{r}\down{l}} 
    \tensor{S}{\down{j}\down{k}} & - \\
  \frac{-1 + 2\,\nu}{24\,\left( -1 + \nu \right) \,\left( 1 + \nu \right) } 
   \,\,\tensor{U}{\down{i}\down{k}\down{r}} 
    \tensor{S}{\down{j}\down{l}} & + &
  \frac{\nu\,\left( -1 + 2\,\nu \right) }
    {24\,\left( -1 + \nu \right) \,\left( 1 + \nu \right) } 
   \,\,\tensor{U}{\down{i}\down{r}\down{k}} 
    \tensor{S}{\down{j}\down{l}} & - \\
  \frac{-1 + 2\,\nu}{24\,\left( -1 + \nu \right) } 
   \,\,\tensor{U}{\down{i}\down{k}\down{l}} 
    \tensor{S}{\down{j}\down{r}} & + &
  \frac{-1 + 2\,\nu}{24\,\left( -1 + \nu \right) } 
   \,\,\tensor{U}{\down{i}\down{l}\down{k}} 
    \tensor{S}{\down{j}\down{r}} & - \\
  \frac{{\left( -1 + 2\,\nu \right) }^2}
    {24\,\left( -1 + \nu \right) \,\left( 1 + \nu \right) } 
   \,\,\tensor{U}{\down{i}\down{j}\down{l}} 
    \tensor{S}{\down{k}\down{r}} & - &
  \frac{\left( -2 + \nu \right) \,\left( -1 + 2\,\nu \right) }
    {24\,\left( -1 + \nu \right) \,\left( 1 + \nu \right) } 
   \,\,\tensor{U}{\down{i}\down{l}\down{j}} 
    \tensor{S}{\down{k}\down{r}} & + \\ 
  \frac{{\left( -1 + 2\,\nu \right) }^2}
    {24\,\left( -1 + \nu \right) \,\left( 1 + \nu \right) } 
   \,\,\tensor{U}{\down{i}\down{j}\down{k}} 
    \tensor{S}{\down{l}\down{r}} & + &
  \frac{\left( -2 + \nu \right) \,\left( -1 + 2\,\nu \right) }
    {24\,\left( -1 + \nu \right) \,\left( 1 + \nu \right) } 
   \,\,\tensor{U}{\down{i}\down{k}\down{j}} 
    \tensor{S}{\down{l}\down{r}} & \\
\end{array}
\]
This sum has the structure
\begin{eqnarray} \label{equ5.8}%
\mathfrak{P}_{i_1\ldots i_5}^{\rm red} & = & \sum_{q\in\calS_5}\,\frac{P_q^{\calP}(\nu)}{Q_q^{\calP}(\nu)}\,U_{i_{q(1)} i_{q(2)} i_{q(3)}} S_{i_{q(4)} i_{q(5)}}\,.
\end{eqnarray}
where $P_q^{\calP}(\nu)$ and $Q_q^{\calP}(\nu)$ are polynomials. 
If we determine the set $N_{\calP}$ of all roots $\nu\not= \frac{1}{2}$ of the $P_q^{\calP}(\nu)$, for which ${\Delta}_{\calP}(\nu)\not= 0$, and set the $\nu\in N_{\calP}$ into (\ref{equ5.8}), the length of (\ref{equ5.8}) will decrease.

We determine the minimal length of (\ref{equ5.8}) by this procedure for $y_t^{\ast}(U\otimes w)$, $y_{t'}^{\ast}(U\otimes S)$, $y_{t'}^{\ast}(U\otimes A)$ and for every of the 15 identity systems of type (\ref{equ5.7}) in the case $\nu\not= \infty$. Table 1 shows the results for $y_{t'}^{\ast}(U\otimes S)$, $y_{t'}^{\ast}(U\otimes A)$. (In the column for $\calP = \{p_1,p_2\}$ the $p_i$ are denoted by their numbers in the lexicographically ordered $\calS_3$.) Furthermore we calculate the lengths of (\ref{equ5.8}) for the 15 systems (\ref{equ5.7}) in the case $\nu = \infty$. Altogether, the number of calculations comes to
\centerline{(3 generator types)$\times$ (2 $\nu$-cases)$\times$ (15 systems (\ref{equ5.7})) = 75.}
We obtain (see \cite{fie03c,fie04a})

\begin{theorem}{Theorem} \label{thm12}%
Let $\dim V\ge 3$. Then the coordinates of $y_t^{\ast}(U\otimes w)$, $y_{t'}^{\ast}(U\otimes S)$, $y_{t'}^{\ast}(U\otimes A)$ are sums of the following lengths
\begin{center}
\begin{tabular}{|c|l|c|c|c|}
\hline
 & & $y_t^{\ast}(U\otimes w)$ & $y_{t'}^{\ast}(U\otimes S)$ & $y_{t'}^{\ast}(U\otimes A)$ \\
\hline
{\rm (a)} & {\rm generic case for $\nu$} & {\rm 8} & {\rm 16} & {\rm 20} \\
\hline
{\rm (b)} & {\rm $\nu$ producing minimal length} & {\rm 4} & {\rm 12} & {\rm 10} \\
\hline
\end{tabular}
\end{center}
\end{theorem}
The computer calculations were carried out by means of the Mathematica packages {\sf Ricci} \cite{ricci3} and {\sf PERMS} \cite{fie10}. Notebooks of the calculations are available on the web page \cite{fie21}.

It is very remarkable that $U$ admits an index commutation symmetry if the coordinates of $y_t^{\ast}(U\otimes w)$, $y_{t'}^{\ast}(U\otimes S)$, $y_{t'}^{\ast}(U\otimes A)$ have the minimal lengths of case {\rm (b)} in Theorem \ref{thm12} (see \cite{fie03c,fie04a}).

{\small
\begin{center}
\begin{tabular}{c|c|c|c|c}
\hline
 & \multicolumn{2}{c}{$y_{t'}^{\ast}(U\otimes S)$}\vline & \multicolumn{2}{c}{$y_{t'}^{\ast}(U\otimes A)$}\\
$\calP$ & 
roots of $P_q^{\calP}(\nu)$ with & 
length of  & 
roots of $P_q^{\calP}(\nu)$ with  & 
length of \\
 & 
${\Delta}_{\calP}(\nu)\not= 0\;,\;\nu\not= 1/2$ & 
$\mathfrak{P}_{i_1\ldots i_5}^{\rm red}$ &  
${\Delta}_{\calP}(\nu)\not= 0\;,\;\nu\not= 1/2$ & 
$\mathfrak{P}_{i_1\ldots i_5}^{\rm red}$ \\
\hline
\hline
 12 &  0 & 12 &    & \\
    &  2 & 14 & 2  & 12 \\
\hline
 13 & -1 & 14 & -1 & 18 \\
\hline
 14 & -1 & 14 & -1 & 18 \\
    &  0 & 12 &    & \\
    &  2 & 12 &  2 & 10 \\
\hline
 15 & -1 & 12 & -1 & 10 \\
    &  0 & 12 &    & \\
    &  1 & 12 &    & \\
    &  2 & 14 &  2 & 12 \\
\hline
 16 & -1 & 12 & -1 & 10 \\
    &  0 & 12 &    & \\
    &  1 & 12 &    & \\
    &  2 & 12 &  2 & 10 \\
\hline
 23 & -1 & 14 & -1 & 18 \\
    &  0 & 12 &    & \\
    &  2 & 14 &  2 & 12 \\
\hline
 24 & -1 & 14 & -1 & 18 \\
    &  2 & 14 &  2 & 18 \\
\hline
 25 & -1 & 12 & -1 & 10 \\
    &  1 & 12 &    & \\
\hline
 26 & -1 & 12 & -1 & 10 \\
    &  1 & 12 &    & \\
    &  2 & 14 &  2 & 18 \\
\hline
 34 &  0 & 12 &    & \\
    &  2 & 12 &  2 & 10 \\
\hline
 35 & -1 & 14 & -1 & 12 \\
    &  0 & 12 &    & \\
    &  1 & 12 &    & \\
    &  2 & 14 &  2 & 12 \\
\hline
 36 & -1 & 14 & -1 & 12 \\
    &  0 & 12 &    & \\
    &  1 & 12 &    & \\
    &  2 & 12 &  2 & 10 \\
\hline
 45 & -1 & 14 & -1 & 12 \\
    &  1 & 12 &    & \\
    &  2 & 14 &  2 & 18 \\
\hline
 46 & -1 & 14 & -1 & 12 \\
    &  1 & 12 & \\
\hline
 56 &  2 & 14 &  2 & 18 \\
\hline
\hline
\end{tabular}\\
\vspace{3mm}
\uppercase{Table} 1. The lengths of $\mathfrak{P}_{i_1\ldots i_5}^{\rm red}$ for $y_{t'}^{\ast}(U\otimes S)$, $y_{t'}^{\ast}(U\otimes A)$ and $\nu\not= \infty$.
\end{center}
}


\end{document}